\documentclass[a4paper,reqno,10pt]{amsart}
\usepackage{t1enc} 
 
\usepackage[latin1]{inputenc}
\usepackage[english]{babel}
\usepackage{amsfonts} 
\usepackage{enumerate}
\usepackage{epsfig}
\usepackage{graphicx}
\usepackage{psfrag}
\usepackage{amsmath}
\usepackage{amsthm}
\usepackage{amssymb}
\usepackage{amscd} 
\usepackage{stmaryrd}

\usepackage{bbm} 
\usepackage{mathrsfs} 

\usepackage{verbatim} 

\newtheoremstyle{exercise} 
  {3pt} 
  {3pt} 
  {\scriptsize\rmfamily} 
  {
\parindent} 
  {\rmfamily\scshape} 
  {.} 
  {.5em} 
  {} 

\newtheoremstyle{newplain}
  {5pt}
  {5pt}
  {\itshape}
  {}
  {\rmfamily\scshape}
  {. ---}
  {.5em}
  {}

\newtheoremstyle{newremark}
  {5pt}
  {5pt}
  {\rmfamily}
  {}
  {\rmfamily\scshape}
  {. ---}
  {.5em}
  {}

\theoremstyle{newplain}
\newtheorem*{Theorem*}{Theorem} 

\theoremstyle{newplain}
\newtheorem{Theorem}{Theorem}
\newtheorem{Lemma}[Theorem]{Lemma}
\newtheorem{Corollary}[Theorem]{Corollary}
\newtheorem{Proposition}[Theorem]{Proposition}

\newtheorem{Definition}[Theorem]{Definition}

\newtheorem{ques}[Theorem]{Question}

\theoremstyle{newremark}
\newtheorem{Remark}[Theorem]{Remark}

\theoremstyle{exercise}

\numberwithin{Theorem}{section}
\numberwithin{Exercise}{section}

\def\e{\varepsilon}

\def\F1{{\mathscr F}}

\def\Gje{\Gamma_j^{\eta}}

\def\pf{\begin{proof}[{\bf Proof:}]\thst}
\def\thst{\mbox{}\newline}

\def\dimuB{\underline{dim}_B}
\def\dimoB{\overline{dim}_B}
\def\dimhm{dim_{\mathscr{H}}}

\def\dimP{dim_{\mathscr{P}}}

\newcommand{\N}{\mathbb{N}} 
\newcommand{\R}{\mathbb{R}} 

\newcommand{\Rn}{\R^n}
\newcommand{\Hm}[1]{\mathscr{H}^{#1}} 
\newcommand{\Pm}[1]{\mathscr{P}^{#1}} 

\newcommand{\cl}[1]{\mathcal{#1}}

\newcommand{\calA}{\mathscr{A}}

\def\XXint#1#2#3{{%
\setbox0=\hbox{$#1{#2#3}{\int}$}
\vcenter{\hbox{$#2#3$}}\kern-.5\wd0}}

\renewcommand{\leq}{\leqslant}
\renewcommand{\geq}{\geqslant}

\begin{document}


\title{Minkowski and packing Dimension comparisons for sets with Reifenberg properties}
\author{Amos N. Koeller}
\address {Mathematisches Institut der Universit\"at T\"ubingen \\ 
Auf der Morgenstelle 10 \\
72076 T\"ubingen \\
Germany }
\email{akoeller@everest.mathematik.uni-tuebingen.de}
\date{\today}


\begin{abstract}
In Koeller \cite{koerprops} the twelve variants of the Reifenberg properties known to be instrumental in the theory of minimal surfaces were classified with respect to various Hausdorff measure based measure theoretic properties. The classification lead to the consideration of fine geometric properties and a connection to fractal geometry. The current work develops this connection and extends the classification to consider Minkowski-dimension, packing dimension, measure, and rectifiability, and the equality of packing and Hausdorff measures with interesting results.
\end{abstract}

\maketitle

\section{Introduction}
In 1960 Reifenberg \cite{reif} considered sets $A\subset B_{\rho_0}(x)\subset \R^n$ for some $x\in \R^n$ and $\rho_0>0$ for which, for all $y\in A$ and $\rho \in (0,\rho_0]$, there exists a $j$-dimensional plane $L_{y,\rho}$ such that
$d_{\Hm{}}(A\cap B_{\rho}(y),L_{y,\rho}\cap B_{\rho}(y))<\e \rho$. Such sets have become known as sets satisfying the $j$-dimensional $\e$-Reifenberg property. 

Simon \cite{simon2} used a similar property in his important work showing the rectifiability of a class of minimal surfaces, the most important change being that the set should be $\e$-Reifenberg for each $\e>0$.

Further similar properties have been widely investigated, among others by David and Toro \cite{davidtoro}, David, Kenig and Toro \cite{davkentor}, David, de Pauw and Toro \cite{davdeptor}, and de Pauw and Koeller \cite{depkoe}, considering generalisations of Reifenberg's work, shapes of Reifenberg sets, and a graph version of Reifenberg's property. The applications of the varying properties such as those given by Reifenberg in his original work and Simon in \cite{simon2}, as well as of the generalisations such as given by David, de Pauw and Toro \cite{davdeptor} show that an understanding of the structure and properties in question is important.

It is the geometric measure theoretic characteristics of sets satisfying such properties that we investigate here. Exactly which properties should be investigated is made unclear by the fine differences between the definitions given by Reifenberg and Simon. We therefore consider all of the twelve possible affine approximation methods, as defined below in Definition $\ref{defa}$, arising from combinations of the essential ingredients of the approximations given by Reifenberg and Simon. 

The measure theoretic characteristics of the Reifenberg properties that we are interested in are the dimension, locally finite measure and rectifiability of sets satisfying the Reifenberg properties. Especially given, as we shall see, that some of the $j$-dimensional Reifenberg properties do not even ensure that the set be (Hausdorff) $j$-dimensional it is with these very general measure theoretic characteristics that we must start. In this paper we classify the Reifenberg properties with respect to which of the mentioned characteristics are ensured.

In a previous paper, \cite{koerprops}, we have already completed the classification with respect to Hausdorff dimension and measure as well as rectifiability with respect to the Hausdorff measure. In the fine geometric structure of Euclidean spaces, however, Hausdorff measure is not the only important measure. In this paper we extend the classification to include the classification with respect to packing and Minkowski dimensions. The Reifenberg properties are additionally classified with respect to which properties ensure locally finite packing measure, rectifiability with respect to the packing measure and with respect to which properties ensure that the packing and Hausdorff measures agree.

From the classification we also draw further insight into Minkowski dimension and the interplay between the packing and Hausdorff measures. The fact that Minkowski dimension often does not provide the heuristically expected dimension is emphasised in that even sets with very strong $j$-dimensional affine approximations can have Minkowski dimension exceeding $j$. On the other hand we show that some Reifenberg properties ensure that the packing and Hausdorff measures are identical, even when the Hausdorff measure (and therefore also packing measure) is not even locally finite.

To formally define the variants of Reifenberg's property being considered, we first define, for any $A\subset \R^n$ and any non-negative real number $r$, $A^r:=\{x\in \R^n:d(x,A)\}$ (where $d$ denotes the usual Euclidean distance). Our variants of Reifenberg's property can now be defined as the twelve following properties.

\begin{Definition}\label{defa} \thst
Let $A \subset \R^n$ be an arbitrary set and $j \in \mathbb{N}$; then
\begin{enumerate}[(i)]
\item $A$ has the weak $j$-dimensional $\delta$-approximation property (or \emph{$wj$ property}) for some $0<\delta <1$ if, for all $y \in A$, there is a $\rho_y > 0$ such that for all $\rho \in (0,\rho_y]$ there exists $L_{y,\rho}\in G_y(n,j)$ such that $B_{\rho}(y) \cap A \subset L_{y, \rho}^{\delta \rho}$.
\item $A$ has the weak $j$-dimensional $\delta$-approximation property with local $\rho_y$-uniformity (or \emph{$w\rho j$ property}) for some $0<\delta <1$ if, for all $y \in A$, there is a $\rho_y > 0$ such that for all $\rho \in (0,\rho_y]$ and all $x \in B_{\rho_y}(y) \cap A$, there exists $L_{x,\rho} \in G_x(j,n)$ such that $B_{\rho}(x) \cap A \subset L_{x, \rho}^{\delta \rho}$.
\item The property (i) is said to be $\rho_0$-uniform (referred to as the \emph{$w\rho_0 j$ property}), if $A$ is contained in some ball of radius $\rho_0$ and if, for every $y \in A$ and every $\rho \in (0,\rho_0]$, there exists $L_{y,\rho}\in G_y(j,n)$ such that $B_{\rho}(y) \cap A \subset L_{y, \rho}^{\delta \rho}$.
\item $A$ is said to have the fine weak $j$-dimensional approximation property (or \emph{$w\delta j$ property}) if, for each $\delta>0$, A satisfies (i).
\item $A$ is said to have the fine weak $j$-dimensional approximation property with local $\rho_y$-uniformity (or \emph{$w\rho\delta j$ property}) if, for each $\delta>0$, $A$ satisfies (ii).
\item $A$ is said to have the fine weak $j$-dimensional approximation property with \emph{$\rho_0$-uniformity} (or $w\rho_0\delta j$ property) if, for each $\delta>0$, $A$ satisfies (iii).
\item $A$ is said to have the strong $j$-dimensional $\delta$-approximation property (or \emph{$s j$ property}) for some $0<\delta <1$ if, for each $y \in A$, there exists $L_y \in G_y(j,n)$ such that definition (i) holds with $L_{y, \rho} = L_y$ for every $\rho \in (0,\rho_y]$.
\item $A$ is said to have the strong $j$-dimensional $\delta$-approximation property (or \emph{$s\rho j$ property}) with local $\rho_y$-uniformity for some $0<\delta <1$ if, for all $y \in A$, there exists $L_y\in G(j,n)$ such that for all $x \in B_{\rho_y}(y)$ and all $\rho \in (0,\rho_y]$ we have $B_{\rho}(x) \cap A \subset (L_y+x)^{\delta \rho}.$ 
\item The property in (viii) is said to be $\rho_0$-uniform (referred to as the \emph{$s\rho_0 j$ property}) if $A$ is contained in some ball of radius $\rho_0$ and there exists $L\in G(j,n)$ such that for each $x\in A$ and $\rho\in (0,\rho_0]$
$B_{\rho}(x)\cap A \subset (L+x)^{\delta\rho}.$
\item $A$ is said to have the fine strong $j$-dimensional approximation property (or \emph{$s\delta j$ property}) if, for each $\delta>0$, $A$ satisfies (vii).
\item $A$ is said to have the fine strong $j$-dimensional approximation property with local $\rho_y$-uniformity (or \emph{$s\rho\delta j$ property}) if, for each $\delta>0$, $A$ satisfies (viii).
\item $A$ is said to have the fine strong $j$-dimensional approximation property with $\rho_0$-uniformity (or \emph{$s\rho_0\delta j$ property}) if, for each $\delta>0$, $A$ satisfies (ix).
\end{enumerate}
Such a property as defined above will be referred to in general as a \emph{$j$-dimensional Reifenberg property} or a \emph{Reifenberg property} if the dimension is clear from the context.

For $\alpha\in \{w,s\}$, $\beta\in \{\emptyset, \rho, \rho_0\}$, $\gamma\in \{\delta\}\cup (0,1)=:\Delta$ and $j\leq n$ we write 
$R(\alpha, \beta,\gamma ; j)$ to denote the set of subsets of $\R^n$ satisfying the $\alpha\beta\gamma j$ property if $\gamma=\delta$ and to denote the set of subsets of $\Rn$ satisfying the $\alpha\beta j$ property with respect to $\gamma$ otherwise.
\end{Definition} \noindent
\begin{Remark}
The three essential elements of the definitions are whether the approximation is weak or strong, whether the approximation is $\rho_y$ uniform, $\rho_0$ uniform, or without local uniformity, and whether the approximation is $\delta$-fine or not. With these three elements, the notation $R(\alpha, \beta,\gamma ; j)$ can be seen to be descriptive of how a set is approximated. For further discussion on the not particularly transparent list of definitions, see \cite{koerprops}.

A difference to the original Reifenberg property of note is that the Reifenberg property is two-sided in the sense that the approximating term is $d_{\Hm{}}(A\cap B_{\rho}(y),L_{y,\rho}\cap B_{\rho}(y))<\e \rho$ instead of $A\cap B_{\rho}(y) \subset L_{y,\rho}^{\delta \rho}$. We consider the one-sided version as it allows for more sets to be considered. This may be done, as motivated in \cite{koerprops}, as it makes no difference to the resulting classification.

Note finally that the motivating property considered by Simon in \cite{simon2} is exactly the $w\rho\delta j$ property. The property originally considered by Reifenberg in \cite{reif} can be stated as the two sided version of the $w\rho_0 j$ property.
\end{Remark}

Before continuing, we note the following simple but important inclusion relations between the Reifenberg properties. 
\begin{Proposition}\label{defarels}
Let $j, n\in \N$, $j\leq n$, $\alpha\in \{w,s\}$, $\beta\in \{\emptyset, \rho, \rho_0\}$ and $\gamma\in \Delta$.
Then
$$R(s, \beta, \gamma; j) \subset R(w, \beta, \gamma; j),$$
$$R(\alpha, \rho_0, \gamma; j) \subset R(\alpha, \rho, \gamma; j)\subset R(\alpha, \emptyset, \gamma; j),$$
$$R(\alpha, \beta, \delta; j) \subset R(\alpha, \beta, \gamma_1; j) \subset R(\alpha, \beta, \gamma_2; j) \hbox{ for }0<\gamma_1\leq \gamma_2\leq 1,  \hbox{ and}$$
$$R(\alpha, \beta, \gamma; j) \subset R(\alpha, \beta, \gamma; j+1).$$
Furthermore, if $A\subset B\in  R(\alpha, \beta, \gamma; j)$, then $A\in R(\alpha, \beta, \gamma; j)$.
\end{Proposition}

\section{Preliminaries}
Having introduced the intention of this paper generally, we now make our setting more precise. We start by establishing our notation.

We write $\Hm{j}$ to denote $j$-dimensional Hausdorff measure, $\omega_j$ to denote the $\Hm{j}$ measure of the $j$-dimensional unit ball, and $\dimhm$ to denote Hausdorff dimension. We write $\dimuB$ and $\dimoB$ to denote the upper and lower Minkowski dimensions respectively and $\dimP$ to denote the packing dimension. The packing measure will be denoted by $\Pm{j}$. Moreover $G(n, j)$ will denote the Grassmann manifold of $j$-dimensional subspaces of $\Rn$.

For full definitions and the basic properties of the above measures and dimensions we refer to Falconer \cite{falconer1} and  \cite{falconer} or Mattila \cite{mattila}. Note, however, that, unlike some works on fractal geometry we use $\omega_j$ to ensure that $\Hm{j}$ corresponds to Lebesgue measure on $\R^j$.

As mentioned, it is with respect to measure theoretic characteristics pertaining to the Minkowski dimension as well as the Packing dimension and measure that we give our classification. Apart from the Minkowski and Packing dimension of sets approximated by a given Reifenberg property, the characteristics with which we concern ourselves are defined below.

\begin{Definition}\label{locfinandrect}
Let $j,n\in \mathbb{N} \cup\{0\}$, $j\leq n$, $\mu\in \{\Hm{j}, \Pm{j}\}$ and $A\subset \R^n$. 

$A$ will be said to have strongly locally finite $\mu$ measure (or strong local $\mu$-finality) if for all compact subsets $K \subset \R^n$, $\mu(K \cap A) <\infty ,$ or equivalently, if for all $y \in \R^n$ there exists a radius $\rho_y > 0$ such that $\mu(B_{\rho_y}(y) \cap A) <\infty.$

$A$ will be said to have weakly locally finite $\mu$ measure (or weak local $\mu$-finality) if, for each $y \in A$, there exists a radius $\rho_y > 0$ such that $\mu(B_{\rho_y}(y) \cap A ) < \infty.$

$A$ is said to be ($\mu$, $j$)-rectifiable if and only if 
\[A\subset M_0 \cup \bigcup_{i=1}^{\infty}f_i(\R^j)\]
where $\mu(M_0)=0$ and $f_i$ is a Lipschitz function for each $i \in \N$. Finally, $A$ is said to be purely ($\mu$, $j$)-unrectifiable if for all ($\mu$,$j$)-rectifiable subsets, $F\subset A$, $\mu(F)=0$.
\end{Definition}
\begin{Remark}
Both definitions of locally finite measure are considered in the literature, and therefore both are used here in what is intended to be a general classification. Further comment on the differing definitions and their importance for Reifenberg properties can be found in \cite{koerprops}. Our definition of rectifiability follows essentially the definition of $(\mu,j)$-countably rectifiable sets of Federer \cite{federer}.
\end{Remark}
We are now in a position to formulate formally the questions by which we classify the Reifenberg properties. 
\begin{ques}\label{quesnew}
For each $\alpha\in \{w,s\}$, $\beta\in \{\emptyset, \rho, \rho_0\}$ and $\gamma\in \Delta$, does $P\in R(\alpha,\beta,\gamma;j)$ imply that $P$
\begin{enumerate}
\item has lower (respectively upper) Minkowski dimension less than or equal to $j$,
\item has packing dimension less than or equal to $j$,
\item has - (a) weakly or (b) strongly - locally finite $\Pm{j}$-measure,
\item is ($\Pm{j}$, $j$)-rectifiable,
\item satisfies $\Pm{j}|_P=\Hm{j}|_P$ on $\Rn$
\end{enumerate}
for each $j\in \N$?
\end{ques} \noindent
\begin{Remark}
Our classification will answer these questions. We formulate such an answer by saying that the answer to a given property and question is either yes or no. For example, the $j$ dimensional property (i) allows for sets satisfying (i) with packing dimension greater than $j$. We therefore say that the answer to (i) (1) is no.
\end{Remark}
As has been mentioned, the analogous questions to (2), (3), and (4) for Hausdorff measure have been answered in \cite{koerprops}. The classification for Hausdorff measure becomes a useful tool of reference in this work. The questions for the Hausdorff classification can be stated as below.
\begin{ques}\label{quesold}
For each $\alpha\in \{w,s\}$, $\beta\in \{\emptyset, \rho, \rho_0\}$ and $\gamma\in \Delta$, does $P\in R(\alpha,\beta,\gamma;j)$ imply that $P$
\begin{enumerate}
\item has Hausdorff dimension less than or equal to $j$?,
\item has - (a) weakly or (b) strongly - locally finite $\Hm{j}$-measure?,
\item is ($\Hm{j}$, $j$)-rectifiable?
\end{enumerate}
for each $j\in \N$?
\end{ques} \noindent
For reference we provide the classification of the Reifenberg properties with respect to Question $\ref{quesold}$
\begin{Theorem}\label{hausclassification}
The properties defined in Definition $\ref{defa}$ satisfy the classification with respect to the questions given in Question $\ref{quesold}$ given in the table below.
\begin{equation} \label{haustab}
\begin{tabular}[h]{lccc}

\hbox{Property} &  & \hbox{Question} &     \nonumber \\
\hline
& & &  \nonumber \\
& (1) & (2) & (3) \nonumber \\
&  &  \hbox{ (a), (b)}&     \nonumber \\ \hline
& & &  \nonumber \\
$w j$ & \hbox{No} & \hbox{No, No}  & No  \nonumber \\
$w\rho j$  & \hbox{No} & \hbox{No, No} & No  \nonumber \\
$w\rho_0 j$ & \hbox{No}  & \hbox{No, No} & No  \nonumber \\
$w\delta j$ & \hbox{Yes} & \hbox{No, No}& No  \nonumber \\
$w\rho\delta j$ & \hbox{Yes} &  \hbox{No, No}& No  \nonumber \\
$w\rho_0\delta j$ & \hbox{Yes} & \hbox{Yes, Yes} & Yes  \nonumber \\
$s j$ & \hbox{Yes}  & \hbox{No, No} & Yes \nonumber \\
$s\rho j$ & \hbox{Yes} & \hbox{Yes, No}& Yes  \nonumber \\
$s\rho_0 j$ & \hbox{Yes} &  \hbox{Yes, Yes}& Yes  \nonumber \\
$s \delta j$ & \hbox{Yes} & \hbox{No, No} & Yes  \nonumber \\
$s \rho\delta j$& \hbox{Yes}  & \hbox{Yes, No} & Yes \nonumber \\
$s\rho_0\delta j$ & \hbox{Yes} & \hbox{Yes, Yes} & Yes  \nonumber \\ \hline

\end{tabular}
\end{equation}
\end{Theorem}

Of particular interest in \cite{koerprops} in the analysis leading to the above classification is the following rectifiability theorem, which is integral both to the above classification and the current work.
\begin{Proposition}\label{lipgraphrep}
Let $j \leq n \in \N$. Suppose $A\subset \R^n$ satisfies either $A\in R(w,\rho_0,\delta;j)$ or a strong Reifenberg property, that is $A\in R(s,\beta,\gamma ;j)$ for some $\beta$ and $\gamma$. Then 
\begin{equation}
A\subset \bigcup_{k=1}^{\infty}G_k
\nonumber
\end{equation}
where $G_k$ is the graph of some Lipschitz function over some $j$-dimensional plane. Furthermore, the representation on the right hand side can be taken to be a finite union if $A$ additionally satisfies a $\rho_0$-uniform property. That is, if $A\in R(w,\rho_0,\delta;j)$ or $A\in R(s,\rho_0,\gamma ;j)$ for some $\gamma$.
\end{Proposition}
\begin{proof}
The result follows from 3.9 Lemma 1 in Simon \cite{simon3} and Corollary 3.2 in Koeller \cite{koerprops}.
\end{proof}
Also relevant to the present work from the analysis in \cite{koerprops} is the construction of a particular family of counter examples which assists in showing the irregular results. That is, in showing that the answers to some particular questions is no. From the analysis of this family of sets we deduce the following consequence:
\begin{Proposition}\label{ganda}
Let $j \in \N$ and $0<\eta<1$. Then there exist sets $\Gje \in R(w,\rho_0,\eta;j)$ and $\calA_j \in R(w,\rho,\delta;j)$, each satisfying the two-sided version of the respective Reifenberg properties ((iii) and (v)) such that
\begin{enumerate}[(i)]
\item $\Gje , \calA_j \subset \R^n$ for each $j+1\leq n\in \N$,
\item $\dimhm \Gje > j$,
\item $\calA_j$ has neither weakly nor strongly locally finite $\Hm{j}$-measure, and 
\item $\calA_j$ is purely ($\Hm{j}$, $j$)-unrectifiable.
\end{enumerate} 
\end{Proposition}
\begin{Remark}
The construction of these sets, as well as proofs of their relevant properties can be found in \cite{koerprops}.

The set $\Gje$ can be taken to be the well known fractal, the Koch curve (with appropriate initial angles). Although again related to the Koch curve, $\calA_j$ allows more variability in its construction and, in particular more flatness at fine scales.
\end{Remark}

To conclude this section we recall a couple of simple facts about packing measures and dimension that we find direct use for in this work. Proofs of these results can be found, for e.g., in Mattila \cite{mattila}

\begin{Proposition}\label{packdecomp}
$$dim_{\Pm{}}A=\overline{dim}_PA:=\inf\left\{\sup_i \overline{dim}_MA_i:A=\bigcup_{i=1}^{\infty}A_i,A_i \hbox{ is bounded}\right\}.$$
Furthermore, for $s\in\R$, $n\in \mathbb{N}$ and $s\leq n$
\begin{enumerate}
\item $\Hm{s}\leq \Pm{s}$ on $\R^n$ and
\item should $E \subset \R^n$ satisfy $0<\Pm{s}(E)<\infty$, then $\Pm{s}(E)=\Hm{s}(E)$ if and only if $s$ is an integer, $\Pm{s}|_{E} << \Hm{s}$ and $E$ is ($\Pm{s}$, $s$)-rectifiable.
\end{enumerate}
\end{Proposition}

\section{Minkowski dimension}
We will see that we can only appropriately bound the Minkowski dimension of a Reifenberg-like set when the set is strongly controlled by Lipschitz functions. This fact will also be later exploited. Based on the known Lipschitz representations of sets satisfying certain Reifenberg properties, we are actually able to prove the result giving the desired classification with respect to Minkowski dimension directly.

\begin{Theorem}\label{minkowskidim}
Suppose $j,n\in \N$, $j \leq n$ and $A\subset \R^n$. 

Should $A\in R(w,\rho_0,\delta;j)$, $R(s,\rho_0,\delta;j)$ or $R(s,\rho_0,\gamma ;j)$ for some $\gamma$, then 
$$\dimuB A \leq \dimoB A \leq j.$$
Conversely, if $j<n$, there exists a set, $\cl{N}_j$, satisfying each of the $9$ remaining $j$-dimensional Reifenberg properties for which 
$$\dimoB \cl{N}_j \geq \dimuB \cl{N}_j > j.$$

\end{Theorem}
\begin{proof}
We first consider the case that $A\in R(w,\rho_0,\delta;j)$, $R(s,\rho_0,\delta;j)$ or $R(s,\rho_0,\gamma ;j)$ for some $\gamma$.

We note that for any $j$-dimensional Lipschitz graph over a $j$-dimensional ball, $G$,
\begin{equation}\label{graphminj}
\dimuB G= \dimoB G =j.
\end{equation}
This follows from standard theory. See for e.g. Theorem 5.7 in Mattila \cite{mattila}.

By the definition of the Reifenberg properties, we see that there is a $y\in A$ and a $\rho_0 >0$ such that $A\subset B_{\rho_0}(y)$. Furthermore, by Proposition $\ref{lipgraphrep}$
\begin{equation}\label{mde2}
A\subset \bigcup_{i=1}^Q graph(g_i)
\end{equation}
for Lipschitz functions $g_i:L_i \rightarrow L_i^{\perp}$ and $j$-dimensional subspaces $L_i$ of $\R^n$. As it follows that $A\cap graph(g_i) \subset graph(g_i)|_{B_{\rho_0}(y_i)\cap L_i}$, we deduce from ($\ref{graphminj}$) that 
$$\dimoB(A\cap (graph(g_i))) \leq j$$
for each $i\in \{1,...,Q\}$. Since the union in ($\ref{mde2}$) is finite, we infer that 
$$\dimoB A =\dimoB\left(\bigcup_{i=1}^Q(A\cap(graph(g_i)))\right) \leq j.$$

For the remaining properties, we note that 
$$\cl{N}_j:=\bigcup_{i=1}^{\infty}\{n^{-1}\} \times [0,1]^j \subset \R^n$$
satisfies the $9$ remaining $j$-dimensional Reifenberg properties (for details on the proof, see \cite{koerprops}). We note furthermore, that 
$$\dimuB\cl{N}_0 =2^{-1}$$
(where here $[0,1]^0:=\{0\}$, for a proof, see for e.g. Example 3.5 in Falconer \cite{falconer}). Following now the Notation of Falconer for the definition of Minkowski dimension, that is 
$$\dimuB B:= \liminf_{\varepsilon \rightarrow 0}\frac{ln N_{\varepsilon}(B)}{-ln \varepsilon},$$
we note that for each $0 < \e <<1$
$$N_{\e}(\cl{N}_j) \geq (4\e)^{-j}N_{\e}(\cl{N}_0)$$
so that 
$$\dimoB \geq \dimuB \cl{N}_j \geq \lim_{\e \rightarrow 0}\frac{-jln(4\e) + ln N_{\e}(\cl{N}_0)}{-ln \e} = j+2^{-1}>j.$$
\end{proof}
\begin{Remark}
The above Theorem shows also that the complete classification for either lower or upper Minkowski dimension is identical. We therefore give the classification only once.
\end{Remark}

\section{Packing Dimension}
In the case of packing dimension, the properties for which $dim_{\Pm{}}A>j$ is possible can be found quickly by using known relationships between $dim_{\Pm{}}$ and $dim_{\Hm{}}$ as well as counter examples to bounded dimension known for Hausdorff dimension. For the properties ensuring Lipschitz representation, the results follow similarly quickly from known relationships between $\overline{dim}_M$ and $dim_{\Pm{}}$ as well as ($\ref{graphminj}$). The remaining two properties are the interesting ones. For these two properties, we show that there is a function, $\eta$, dependent on the approximating $\delta$ being used in the Reifenberg property, which yields an upper bound for the dimension of the set. The classification results can then be deduced. We first show the existence of the necessary $\eta$.
\begin{Lemma}\label{betadim}
There exists a function $\eta:(0,\infty) \rightarrow \R$ with $\lim_{\delta\rightarrow 0}\eta(\delta)=0$ such that 
$$dim_{\Pm{}}A\leq j+\eta(\delta_1)$$
whenever $A\in R(w,\rho,\delta_1;j)$.
\end{Lemma}
\begin{proof}
For $\delta_1>1/8$ define $\eta(\delta_1)=n-j+1$. Now, let $0<\delta_1<1/8$ and, for each $x\in A$, define $r_x^{\prime}$ to be the radius for which the $\delta_1$-approximations of the $w\rho j$ property hold around $x$. Define $r_x:=min\{1,r_x^{\prime}\}$. By the Besicovitch covering Theorem there is a countable set $\{x_j\}_{j\in \mathbb{N}}$ such that $A \subset \cup_{j=1}^{\infty}B_{r_{x_j}}(x_j)$. Define $A_j:=B_{r_{x_j}}(x_j) \cap A$.

Note now that for any $j$-dimensional subspace of $\R^n$ and any $\delta> 0$ there exists a constant $C$ such that for each $y\in \R^n$, there exists a cover $\{B_{4\delta}(y_k)\}_{k=1}^Q$ of $L^{2\delta}\cap B_1(y)$ such that $Q(4\delta)^j<C$. Defining for $0<\delta \leq 1/8$
\begin{equation}\label{beta}
\eta(\delta):=-ln(2C)(ln(4\delta))^{-1}
\end{equation}
we see that $\lim_{\delta \rightarrow 0}\eta(\delta)=0$ and $Q(4\delta)^{j+\eta(\delta)}<2^{-1}$. Rescaling we deduce that for each $y\in \R^n$, $0<\delta \leq 1/8$, and $R>0$ there exists a cover $\{B_{4\delta R}(y_k)\}_{k=1}^Q$ of $L^{2\delta R}\cap B_R(y)$ with 
\begin{equation}\label{packeqn1}
Q(4\delta R)^{j+\eta(\delta)}<2^{-1}R^{j+\eta(\delta)}.
\end{equation}
Take some initial cover of $A_j$ with $Q'$ balls of radius $0<\lambda<2^{-1}r_{x_j}$, each with non-empty intersection with $A_j$. Define $T_0:=\max\{Q',Q\}\lambda^{j+\eta(\delta_1)}$.

Now for each of these balls $B_i$ choose $z_i \in B_i \cap A_j$. $B_i \subset B_{2\lambda}(z_i)$ and thus there exists a $j$-dimensional plane $L_i$ with $A_j \cap B_i \subset L^{2\delta_1 \lambda}\cap B_i$ and therefore, from $(\ref{packeqn1})$, there exists a cover of $A_j \cap B_i$ with $Q$ balls of radius $4\delta_1 \lambda$ satisfying 
$$Q(4\delta_1 \lambda)^{j+\eta(\delta_1)}\leq 2^{-1}\lambda^{j+\eta(\delta_1)}.$$
Repeating this argument for each ball $B_i$ we deduce that $A_j$ can be covered by $Q'Q$ balls of radius $4\delta_1 \lambda$ satisfying $Q'Q(4\delta_1 \lambda)^{j+\eta(\delta_1)}\leq 2^{-1}T_0$. That is, 
$$N(A_j,4\delta_1 \lambda)(4\delta_1 \lambda)^{j+\eta(\delta_1)}\leq 2^{-1}T_0.$$

Continuing this process inductively reveals that for each $q\in \mathbb{N}$
$$N(A_j,(4\delta_1)^q\lambda)((4\delta_1)^q\lambda)^{j+\eta(\delta_1)}\leq 2^{-q}T_0.$$
Now, for any $\e>0$ there exists $q\in \mathbb{N}$ such that $(4\delta_1)^{q+1}\lambda < \e < (4\delta_1)^q\lambda$ and as $N(A_j,\e)\leq N(A_j,(4\delta_1)^{q+1}\lambda)$
$$N(A_j,\e)\e^{j+\eta(\delta_1)}\leq N(A_j,(4\delta_1)^{q+1}\lambda)(4\delta_1)^{q+1}\lambda(4\delta_1)^{-1} \leq (2^{q+1}4\delta_1)^{-1}T_0.$$
It follows that $\limsup_{\e \rightarrow 0}N(A_j,\e)\e^{j+\eta(\delta_1)}=0$, hence $\overline{dim_M}A_j\leq j+\eta(\delta_1)$ and therefore that 
\begin{equation}\label{deltavpack}
dim_{\Pm{}}A \leq \sup_{j\in \mathbb{N}}\overline{dim_M}A_j \leq j+\eta(\delta_1).
\end{equation}
Here, the first inequality is a standard characterisation of packing dimension. A proof can be found, for e.g., in Mattila \cite{mattila}.
\end{proof}
Having established the function $\eta$, we are now able to prove the results necessary to classify the Reifenberg properties with respect to packing dimension.
\begin{Theorem}\label{packingtheorem}
Let $A\subset \R^n$ and $j\leq n$. If $A\in R(w,\beta,\delta;j)$ or $A\in R(s,\beta,\gamma;j)$ for some $\beta$ and $\gamma$, then $\dimP A \leq j$. Otherwise $\dimP A>j$ is possible whenever $j<n$. 
\end{Theorem}
\begin{proof}
Since $dim_{\Hm{}}A\leq dim_{\Pm{}}A$ for all $A\subset \R^n$, $dim_{\Hm{}}\Gje >j$ for each $0<\eta<1$ and $j<n$, and $\Gje \in R(w,\rho_0,\eta;j)$, the claim for the $wj$, $w\rho j$, and $w\rho_0 j$ properties follows from Proposition $\ref{defarels}$.

Since, by Proposition $\ref{lipgraphrep}$, any set, $A\subset \R^n$, $A\in R(w,\rho_0,\delta;j)$ or $A\in R(s,\beta,\gamma;j)$ for some $\beta$ and $\gamma$, can be represented by 
$$A\subset \bigcup_{k=1}^{\infty}G_k$$
where $G_k$ is the graph of some Lipschitz function over some $j$-dimensional plane, it follows from Proposition 
$\ref{packdecomp}$ that
\begin{equation}\label{packpropeqn1}
\dim_{\Pm{}}A \leq \sup_{k\in \mathbb{N}}dim_{\Pm{}}G_k.
\end{equation}
By Theorem $\ref{minkowskidim}$, that is, by ($\ref{graphminj}$), we see that $\overline{dim}_MG_k \leq j$ for each $k\in \mathbb{N}$. Thus, since, for any $B\subset \R^n$, $dim_{\Pm{}}B\leq \overline{dim}_M B$, we deduce from ($\ref{packpropeqn1}$) that
$$\dim_{\Pm{}}A \leq \sup_{k\in \mathbb{N}}\overline{dim}_M G_k \leq j.$$
Suppose now that $A\in R(w,\rho,\delta;j)$. Then $A\in R(w,\rho,\e;j)$ for each $\e>0$, and it follows from Lemma $\ref{betadim}$ that 
$$dim_{\Pm{}}A\leq \lim_{\e \rightarrow 0}j+\eta(\e)=j,$$
proving the result claimed for the $w\rho j$ property.

Finally, suppose $A\in R(w,\emptyset,\delta;j)$. Let $s>0$ and take $\e>0$ such that $\eta(\e)< s$ where $\eta:\R\rightarrow \R$ is the function given in ($\ref{beta}$). Since $A\in R(w,\emptyset,\delta;j)$ we know that 
$$\rho_{\e ,x}:=\frac{1}{2}\sup \{r\in \R:r\in R_x\} >0$$
where, for each $x\in A$, $R_x$ is the set of real numbers $\rho_0>0$ such that for all $\rho\in(0,\rho_0]$ there exists a $j$-dimensional affine plane $L_{x,\rho}$ such that $B_{\rho}(x)\cap A \subset L_{x,\rho}^{\e\rho}$.

Define now, for each $m\in \mathbb{N}$, $A_{m}:=\{x\in A:\rho_{\e , x}\geq m^{-1}\}$. Clearly $A=\cup_{m\in \mathbb{N}}A_{m}$. Further, for any $m\in \mathbb{N}$, since $A_{m} \subset A$ and $\rho_{\e,x}$ is bounded below in $A_{m}$, we see that $A_{m}\in R(w,\rho,\e;j)$ with $\rho_y \geq \frac{1}{m}$ for each $y\in A_{m}$. It follows from Lemma $\ref{betadim}$ that $dim_{\Pm{}}A_m \leq j+\eta (\e)$, and thus that
$\Pm{j+s}(A_{m})=0$. Since $m$ was arbitrary
$$0\leq \Pm{j+s}(A)\leq\sum_{m\in \mathbb{N}}\Pm{j+s}(A_m)=0$$ 
and we infer that $dim_{\Pm{}}(A)=\inf\{s\in R:\Pm{s}(A)=0\}\leq j$ which completes the proof.
\end{proof}

By comparing Theorem $\ref{packingtheorem}$ to Theorem $\ref{hausclassification}$ we see that the classifications for Hausdorff and packing dimensions are identical. The classification gives an upper bound on the dimension of the sets, but not a value. It could therefore be asked whether we may conclude that $dim_{\Hm{}}A=dim_{\Pm{}}A$ for a set $A$ satisfying some $j$-dimensional Reifenberg property. 

When relying only on the Reifenberg properties, this is not possible. Even in the case shown in the next section where $\Hm{j}|_A=\Pm{j}|_A$, we cannot say anything about the dimension if the $j$-dimensional Hausdorff and packing measures are zero. To see that the dimensions need not agree, not even for the two sided properties, take first a subset, $A \subset\R^j$ with $0<dim_{\Hm{}}A<dim_{\Pm{}}A$, whose existence is known (see, for e.g. Tricot \cite{tricot}), and define $B:=A\cup \mathbb{Q}^j$. It is clear that $B$ possesses all two-sided $j$-dimensional Reifenberg properties and that $0<\dim_{\Hm{}}B = dim_{\Hm{}}A < dim_{\Pm{}}A=dim_{\Pm{}}B.$ Note, though, that the $j$-dimensional measures need not disagree, as they may both be zero.

Allowing other additional properties, however, can force the Hausdorff and packing dimension to agree. One of the simplest additional properties, as Reifenberg himself considered, is to require the sets to be closed. In this case we can give the following Corollary showing equality of dimension in some cases.

\begin{Corollary}\label{samedimension}
Let $0<\e < \e_n$ where $\e_n$ is the approximating constant given in Reifenberg's topological disc theorem.
Suppose that $A\subset \R^n$ and that
\begin{enumerate}
\item $A\not=\emptyset$,
\item $A$ is closed, 
\item $A\in R(w,\rho_0,\delta;j)$ or $A\in R(s,\beta,\gamma;j)$, and 
\item $A$ satisfies the two sided $w\rho j$ property with respect to $\e$.
\end{enumerate}
Then $dim_{\Pm{}}A=dim_{\Hm{}}A=j$.
\end{Corollary}
\begin{Remark}
We first give more precise definitions of the two-sided properties relevant to this result. 

$A$ satisfies the two-sided $w\rho j$ property with respect to $\e>0$ if for each $y \in A$ there is a $\rho_y>0$ such that for each $x\in A\cap B_{\rho_y}(y)$ and $\rho \in (0,\rho_y]$ there is a $j$-dimensional affine plane $L_{x,\rho}$ satisfying 
$$d_{\Hm{}}(A\cap B_{\rho}(x), L_{x,\rho}\cap B_{\rho}(x))<\e\rho.$$

If $A \subset B_{2\rho_0}(x)$ and $x\in A$, then $A$ satisfies the two sided $w\rho_0 j$ property in $B_{\rho_0}(x)$ if for each $y\in A\cap B_{\rho_0}(x)$ and $\rho \in (0,\rho_0]$ there is a $j$-dimensional affine plane $L_{y,\rho}$ containing $y$ such that 
$$d_{\Hm{}}(A\cap B_{\rho}(y)), L_{y,\rho} \cap B_{\rho}(y))<\e \rho.$$
\end{Remark}
\begin{proof}
Take $y\in A$ and note that $\overline{A\cap B_{\rho_y}(y)}$ is closed and that $A\cap B_{2\rho_y}(y)$ satisfies the two sided version of the $w\rho_0 j$ property in $B_{\rho_y}(y))$, that is, Reifenberg's original condition. By Reifenberg's topological disc theorem it follows that $\overline{A\cap B_{\rho_y}(y)}$ is homeomorphic to a closed $j$-dimensional unit ball (the $j$-dimensional disc).

Since further, by satisfying the two-sided $w\rho j$ property, $d(\partial B_{\rho_y}(y), A)<\e \rho$, we deduce that 
$$\pi_L(A) \supset B_{(1-\e)\rho_y}(y) \cap L_{y,\rho_y},$$
where $\pi_L:\R^n\rightarrow L_{y,\rho_y}$ denotes the projection of $\R^n$ onto $L_{y,\rho_y}$. It follows that 
$$((1-\e)\rho_y)^j\omega_j \leq \Hm{j}(A)\leq \Pm{j}(A)$$
and therefore that $j\leq dim_{\Hm{}}\leq dim_{\Pm{}}$. The result now follows from Theorem $\ref{packingtheorem}$.
\end{proof}

\section{Packing Measure}
Noting again that the classification for Hausdorff dimension and packing dimension are identical adds further to the interest in the packing measure in that we may find nice conditions showing that the packing and Hausdorff measures agree. In fact, we do; we show that in the cases for which there is no appropriate upper bound on the dimension, nothing further can be said about the measure, but that in the remaining cases, the packing and Hausdorff measures, and their classifications agree with each other.

We start by showing the negative results. That is, those questions answered with a no. These results follow quickly from Proposition $\ref{packdecomp}$.
\begin{Lemma}\label{nothnotp}
Let $A\subset \R^n$ and $j \leq n$. Should $A$ not be weakly (respectively strongly) locally $\Hm{j}$-finite, then $A$ is also not weakly (respectively strongly) locally $\Pm{j}$-finite. Similarly, should $A$ not be ($\Hm{j}$, $j$)-rectifiable, then $A$ is also not ($\Pm{j}$, $j$)-rectifiable.
\end{Lemma}
\begin{proof}
The claim regarding locally finite measures follows directly from the definition of locally finite measure and Proposition $\ref{packdecomp}$ (1).

For the second claim, should $A$ not be ($\Hm{j}$, $j$)-rectifiable then for each union 
$$M=\bigcup_{i=1}^{\infty}f_i(\R^j)$$
with Lipschitz functions $f_i$, $\Pm{j}(A\sim M)\geq \Hm{j}(A\sim M)>0$ from which the result follows.
\end{proof}
\begin{Corollary}\label{pnoth}
There exists a set, $A\subset \R^2$ satisfying the $w1$, $w\rho 1$, $w\rho_0 1$, $w\delta 1$, and $w\rho \delta 1$ properties for which
$$\Hm{1}(A)\not= \Pm{1}(A).$$
\end{Corollary}
\begin{proof}
By Proposition $\ref{ganda}$, there is a purely ($\Hm{1}$, $1$)-unrectifiable Borel set $\calA_1 \subset \R^2$, $\calA_1\in R(w,\rho, \delta;1)$, and $\Hm{1}(\calA_1)=\infty$. It follows, see \cite{falconer1} or \cite{falconer}, that there is a compact set $E\subset \calA_1$ satisfying $0<\Hm{1}(E) <\infty$. Since $E\subset \calA_1$, it follows that $E$ is purely ($\Hm{1}$, $1$)-unrectifiable and that $E\in R(w,\rho,\delta;1)$. That $E$ satisfies the $w1$, $w\rho 1$ and $w \delta 1$ properties now follows from Proposition $\ref{defarels}$. 

Now let $\delta >0$ and take, for each $y \in E$, a $\rho_y >0$ such that for each $\rho \in (0,\rho_y]$ and $z\in B_{\rho_y}(y)$ there exists an $L_{z,\rho}\in G(j,n)$ such that
$$E\cap B_{\rho}(z) \subset (L_{z,\rho}+z)^{\delta \rho}.$$
That such a $\rho_y$ exists for each $y\in E$ follows from the definition of the $w\rho \delta 1$ property. Since $E$ is compact we can take $\{y_k\}_{k=1}^Q$ such that 
$$E\subset \bigcup_{k=1}^QB_{\rho_{y_k}}(y_k).$$
We deduce that there exists $y_0 \in \{y_k\}_{k=1}^Q$ such that, defining $A:=E\cap B_{\rho_{y_0}}(y_0)$, $\Hm{1}(A)>0$. It follows that $A$ is purely ($\Hm{1}$, $1$)-unrectifiable. Moreover, by the selection of $\rho_{y_0}$ and that $A\subset E$, we see that $A$ satisfies the $w1$, $w\rho 1$, $w\rho_0 1$, $w\delta 1$, and $w\rho \delta 1$ properties.

Now, if $\Pm{1}(A)=\infty$ the proof is complete. Otherwise, $0<\Hm{1}(A) \leq \Pm{1}(A)<\infty$ and by Lemma $\ref{nothnotp}$ $A$ is not ($\Pm{1}$,$1$)-rectifiable. By Proposition $\ref{packdecomp}$ it follows that $\Pm{1}(A)\not= \Hm{1}(A)$.
\end{proof}
\begin{Remark}
Remarking again on the two-sided situation, we note that the two-sided case also allows for $\Hm{1} \not= \Pm{1}$, indeed consider $A_{\mathbb{Q}}:=\calA_1 \cap \mathbb{Q}_2$ where we define $\mathbb{Q}_2:= ((\R\times \mathbb{Q})\cup(\mathbb{Q}\times \mathbb{R}))$. Since $\mathbb{Q}_2$ is $\Pm{1}$- and ($\Hm{1}$, $1$)-rectifiable, so too is $A_{\mathbb{Q}}$ and we therefore deduce that $\Pm{1}(A_{\mathbb{Q}})=\Hm{1}(A_{\mathbb{Q}})=0$. $A_{\mathbb{Q}}$ is however dense in $\calA_1$ and therefore satisfies the two-sided $w1$, $w\rho 1$, $w\rho_0 1$, $w\delta 1$, and $w\rho \delta 1$ properties. By taking the sets $A$ and $B_{\rho_{y_0}}(y_0)$ found in Corollary $\ref{pnoth}$ and defining $A_2:=A\cup (B_{\rho_{y_0}}(y_0) \cap A_{\mathbb{Q}})$ we see that $A_2$ satisfies the same two sided Reifenberg properties as $A_{\mathbb{Q}}$ but that $\Pm{1}(A_2)\not=\Hm{1}(A_2)$.
\end{Remark}
The positive results, which we now prove, follow from Proposition $\ref{lipgraphrep}$ and from showing that for sets, $A$, satisfying any given one of the properties not addressed in Corollary $\ref{pnoth}$, $\Hm{1}(A)=\Pm{1}(A)$. The equality of Packing and Hausdorff measures follows, as shown below, from known results and some simple estimates on Lipschitz graphs.
\begin{Lemma}\label{hpliplowestimate}
Let $L\in G(j,n)$, $g:L \rightarrow L^{\perp}$ be a Lipschitz function with Lipschitz constant $M$ and $x\in g(L)$. Then there exists a constant $c=c(M,j)$ such that for all $\rho > 0$ 
$$\Hm{j}(B_{\rho}(x)\cap graph(g(L)))\geq c\rho^j.$$
\end{Lemma}
\begin{proof}
Without loss of generality we can assume that $L=\R^j$, $x\in L$, and $g(x)=0$. For $a\in B^j_{\rho(1+M^2)^{-1/2}}(x)$ we calculate $|g(a)|=|g(a)-g(x)|< M\rho(1+M^2)^{-1/2}$. We deduce that $$|(a,g(a))-x|<(\rho^2(1+M^2)^{-1}+M^2\rho^2(1+M^2)^{-1})^{1/2}=\rho$$ and therefore that $$graph\left(g|_{B^j_{\rho(1+M^2)^{-1/2}}(x)}\right)\subset B_{\rho}(x).$$
Writing $G=graph(g|_{B^j_{\rho(1+M^2)^{-1/2}}(x)})$ the result follows since
$$\Hm{j}(G)\geq \Hm{j}\left(B^j_{\rho(1+M^2)^{-1/2}}(x)\right)= \omega_j\rho^j(1+M^2)^{-j/2}=:c(M,j)\rho^j.$$
\end{proof}
\begin{Lemma}\label{hplipupestimate}
Let $L\in G(j,n)$ and $g:L\rightarrow L^{\perp}$ be a Lipschitz function with Lipschitz constant $M$. Let $A\subset graph (g) \subset \R^n$ and $\Hm{j}(A) <\infty$. Then there exists a constant $C=C(M,j)$ such that 
$$\Pm{j}(A) \leq C\Hm{j}(A).$$
\end{Lemma}
\begin{proof}
Let $a:=\Hm{j}(A)$ and $G:=graph (g)$. By rotation we can, without loss of generality, assume that $L=\R^j$. Let $\e>0$ and take $\{B_k\}_{k=1}^Q$, a collection of balls satisfying
$$A\subset \bigcup_{k=1}^QB_k \ \ \hbox{ and } \ \ \sum_{k=1}^Q\omega_j\left(\frac{d(B_k)}{2}\right)^j<a+\e.$$
Now, for any $\eta < \min\{d(B_k):1\leq k \leq Q\}$, consider an $\eta$ packing, $\{C_l\}_{l\in \mathbb{N}}$, of $A$. We see that
$$\bigcup_{l\in \mathbb{N}}C_l \subset \bigcup_{k=1}^Q2B_k,$$
where $2B_k$ is the ball of identical centre to $B_k$ but twice the radius, that the $C_l$ are disjoint and that the $C_l$ are centred on $A$.

Since $g$ is Lipschitz we see that
$$\Hm{j}\left(G\cap\bigcup_{k=1}^Q2B_k\right) \leq \sum_{k=1}^Q\Hm{j}(G\cap 2B_k)\leq 2^jM^j\sum_{k=1}^Q\omega_j\left(\frac{d(B_k)}{2}\right)^j < 2^jM^j(a+\e).$$
Also, by Lemma $\ref{hpliplowestimate}$, $\Hm{j}(G\cap C_l)\geq c(M,j)2^{-j}d(C_l)^j$ for each $l\in \N$. It follows that
\begin{eqnarray}
\sum_{l\in \mathbb{N}}d(C_l)^j  \leq  2^jc(M,j)^{-1}\sum_{l\in \mathbb{N}}\Hm{j}(G\cap C_l) 
\leq  2^jc(M,j)^{-1}\Hm{j}\left(G\cap \bigcup_{k=1}^Q2B_k\right)
\leq  4^jM^jc(M,j)^{-1}(a+\e )\nonumber
\end{eqnarray}
Defining $C:=(4M)^jc(M,j)$, it follows that $P^j_{\eta}(A) \leq C(A+\e )$. By letting $\eta \searrow 0$ and then $\e \searrow 0$ we deduce that $P^j(A) \leq Ca$ and thus that $\Pm{j}(A) \leq Ca$, giving the result.
\end{proof}
\begin{Theorem}\label{pequalsh}
If $A\subset \R^n$ and $A\in R(w,\rho_0,\delta;j)$ or $A\in R(s,\beta,\gamma;j)$ for some $\beta$ and $\gamma$, then $\Pm{j}(A)=\Hm{j}(A).$
\end{Theorem}
\begin{proof}
Since $\Pm{j}$ and $\Hm{j}$ are Borel regular it is sufficient to show that the result holds for Borel sets. We Therefore assume that $A$ is a Borel set. Together with Proposition $\ref{lipgraphrep}$ we may write 
$$A\subset \bigcup_{i=1}^{\infty}A_i$$
where $\{A_i\}_{i\in \N}$ is a pairwise disjoint family of Borel subsets of Lipschitz graphs with $\Hm{j}(A_i) < \infty$. It is clear that $A_i$ is ($\Pm{j}$, $j$)-rectifiable for each $i \in \N$. By Lemma $\ref{hplipupestimate}$ $\Pm{j}(A_i)\leq C(M,j)\Hm{j}(A_i)<\infty$. Moreover, again by Lemma $\ref{hplipupestimate}$, for any $B\subset \R^n$ with $\Hm{j}(B)=0$
$$0\leq \Pm{j}|_{A_i}(B)=\Pm{j}(A_i\cap B)\leq C(M,j)\Hm{j}(A_i\cap B)=0,$$
and thus $\Pm{j}|_{A_i} << \Hm{j}$.

Since $j\in \mathbb{N}$ it follows from Proposition $\ref{packdecomp}$ (2) that $\Pm{j}(A_i)=\Hm{j}(A_i)$ and thus that
$$\Pm{j}(A)=\sum_{i=1}^{\infty}\Pm{j}(A_i)=\sum_{i=1}^{\infty}\Hm{j}(A_i)=\Hm{j}(A).$$
\end{proof}
As a Corollary of Theorem $\ref{pequalsh}$ the remaining classification results may now be shown.
\begin{Corollary}\label{posresults}
Suppose $A\subset \R^n$ and $A\in R(w,\rho_0,\delta;j)$ or $A\in R(s,\beta,\gamma;j)$ for some $\beta$ and $\gamma$. Then $A$ is ($\Pm{j}$,$j$)-rectifiable. 

Furthermore, should $A$ be of weakly (respectively strongly) $\Hm{j}$ locally finite measure, then $A$ is also of weakly (respectively strongly) $\Pm{j}$ locally finite measure.
\end{Corollary}
\begin{proof}
By Proposition $\ref{lipgraphrep}$ it is clear that $A$ is ($\Pm{j}$, $j$)-rectifiable. The results regarding the local finiteness of measure follow directly from Theorem $\ref{pequalsh}$.
\end{proof}
\section{Classification}
All of the necessary results concerning the classification undertaken here have now been proven in the preceding sections. We now collect and state formally our completed classification in the form of a theorem.
\begin{Theorem}\label{packmeasthm}
The properties defined in Definition $\ref{defa}$ can be classified with respect to the questions given in Question $\ref{quesnew}$ as follows:
\begin{equation} \label{tab2}
\begin{tabular}[h]{lccccc}

\hbox{Property} & & & \hbox{Question} &   &  \nonumber \\
\hline
& & & & & \nonumber \\
& (1) & (2) & (3) & (4) & (5) \nonumber \\
& & &  \hbox{ (weak, strong)}& &     \nonumber \\ \hline
& & & & \nonumber \\
(i) & No & \hbox{No} & \hbox{No, No}  & No & No\nonumber \\
(ii) & No & \hbox{No} & \hbox{No, No} & No & No\nonumber \\
(iii) & No & \hbox{No}  & \hbox{No, No} & No & No\nonumber \\
(iv) & Yes & \hbox{Yes} & \hbox{No, No}& No & No \nonumber \\
(v) & No & \hbox{Yes} &  \hbox{No, No}& No & No \nonumber \\
(vi) & No & \hbox{Yes} & \hbox{Yes, Yes} & Yes & Yes \nonumber \\
(vii) & No & \hbox{Yes}  & \hbox{No, No} & Yes & Yes \nonumber \\
(viii) & No & \hbox{Yes} & \hbox{Yes, No}& Yes & Yes \nonumber \\
(ix) & Yes & \hbox{Yes} &  \hbox{Yes, Yes}& Yes & Yes \nonumber \\
(x) & No & \hbox{Yes} & \hbox{No, No} & Yes & Yes \nonumber \\
(xi) & No & \hbox{Yes}  & \hbox{Yes, No} & Yes & Yes  \nonumber \\
(xii) & Yes & \hbox{Yes} & \hbox{Yes, Yes} & Yes & Yes. \nonumber \\ \hline

\end{tabular}
\end{equation}
\end{Theorem}
\begin{proof}
This theorem is a summary of the above results. More specifically the results can be assembled as follows. The answers to question (1) are given in Theorem $\ref{minkowskidim}$. The packing dimension results follow from Theorem $\ref{packingtheorem}$. The results concerning locally finite measure follow from Theorem $\ref{hausclassification}$, Lemma $\ref{nothnotp}$ and Corollary $\ref{posresults}$. The negative rectifiability results follow from Theorem $\ref{hausclassification}$ and Lemma $\ref{nothnotp}$. The positive rectifiability results follow from Corollary $\ref{posresults}$. Finally, the answers to question (4) follow from Corollary $\ref{pnoth}$ and Theorem $\ref{pequalsh}$.
\end{proof}




\begin{thebibliography}{widest-label}
\bibitem[1]{davidtoro} G. David,T. Toro, {  Reifenberg flat metric spaces, snowballs, and embeddings}, Math. Ann. 315 (1999) 641-710.
\bibitem[2]{davkentor} G. David, C. Kenig, T. Toro, { Asymptotically optimally doubling measures and Reifenberg flat sets with vanishing constant}, Comm. Pure Appl. Math. 54 (2001) 385-449.
\bibitem[3]{davdeptor} G. David, T. de Pauw, T. Toro, { A generalization of Reifenberg's theorem in $\R^3$}, Geom. Funct. Anal. 18, No. 4 (2008) 1168-1235
\bibitem[4]{depkoe} T. de Pauw, A. Koeller, { Linearly approximatable functions}, Proc. Amer. Math. Soc.  137  (2009) 1347-1356. 
\bibitem[5]{falconer1} K.J. Falconer, { The geometry of fractal sets}, Cambridge University Press, 1985.
\bibitem[6]{falconer} K.J. Falconer, { Fractal Geometry: Mathematical Foundations and Applications}, John Wiley \& Sons, 1990.
\bibitem[7]{federer} H. Federer, { Geometric Measure Theory}, Springer-Verlag, Berlin-Heidelberg-New York, 1969.
\bibitem[8]{koerprops} A. N. Koeller, { A classification of Reifenberg properties} Submitted 2009
\bibitem[9]{mattila} P. Mattila, { Geometry of Sets and Measures in Euclidean Spaces, Fractals and Rectifiability} Cambridge studies in advanced mathematics 44, Cambridge University Press, 1995.
\bibitem[10]{reif} R.E. Reifenberg, { Solutions of the Plateau problem for $m$-dimensional surfaces of varying topological type} Acta Math. 104 (1960) 1-92.
\bibitem[11]{simon2} L. Simon, { Rectifiability of the Singular Sets of Multiplicity 1 Minimal Surfaces and Energy 
Minimizing Maps.} Surveys in Diff. Geom. {2} (1995) 246-305.
\bibitem[12]{simon3} L. Simon, { Theorems on Regularity and Singularity of Harmonic Maps}, 
ETH Lectures, Birkh\"auser, 1996.
\bibitem[13]{tricot} C. Tricot, { Two definitions of fractional dimension}, Math. Proc. Cambridge Philos. Soc. 91 (1982) 57-74.
\end{thebibliography}
\end{document}